\documentclass[runningheads]{llncs}
\usepackage{amsmath,amsfonts,amssymb,mathtools}
\usepackage{graphicx} 
\usepackage{tikz}
\usepackage{float}
\usetikzlibrary{arrows.meta}
\usepackage{comment}

\newtheorem{observation}{Observation}
\usepackage{url}

\definecolor{orcidgreen}{RGB}{166,206,57}

\begin{document}
\title{Eccentric Connectivity Index of Strongly Connected Digraphs}
\author{Vysakh Chakooth\thanks{Is also a part-time research scholar in the Department of Mathematics, Government College Chittur, Palakkad, Kerala, India - 678104, affiliated to University of Calicut, Thenhipalam, Kerala, India - 673635.}\inst{1}\orcidID{0009-0007-9176-3845}
\and
Prasanth G. Narasimha-Shenoi\inst{2,4}\orcidID{0000-0002-5850-5410}
\and
Prakash G. Narasimha-Shenoi\inst{3,4}\orcidID{0000-0002-6364-8686}
}
\authorrunning{Vysakh Chakooth et al.}
\institute{Department of Mathematics, NSS College Ottappalam,\\
Palakkad, Kerala, India - 679103\\
\email{v2sakh@gmail.com}
\and
Department of Mathematics, Government College Chittur,\\
Palakkad, Kerala, India - 678104\\  \email{prasanthgns@gmail.com}
\and
Department of Mathematics, Maharajas College Ernakulam,\\ Kerala, India - 682011\\
\email{prakashgn@gmail.com}
\and
Department of Collegiate Education, Government of Kerala,\\
Thiruvananthapuram, Kerala India - 695033
}
\maketitle

\begin{abstract}
Let $G = (V, E)$ be a graph with non-empty set of vertices $V$ and set of edges $E$. The \emph{eccentric connectivity index} of the graph $G$ is defined as $$\displaystyle{\xi^C(G) = \sum_{u \in V}  d_u \;ecc(u)}$$ where $d_u$ is the degree and $ecc(u)$ is the eccentricity of the vertex $u \in V$. This article is an attempt to find the \emph{eccentric connectivity index} of  strongly connected digraph $D$ with respect to the metric, \textit{maximum distance} defined by $md(u,v)=\max\{\vec{d}(u,v),\vec{d}(v,u)\}$. An attempt is also made to find the extremal values for strongly connected digraphs. 
\end{abstract}

\section{Introduction}

Topological indices are numerical invariants that represent structural features of graphs and are widely applied in chemical graph theory for predicting physico – chemical and biological properties. They are generally classified into three types: degree–based, distance–based, and distance–cum–degree based indices. Among them, the \emph{eccentric connectivity index}, of simple graphs introduced by Sharma et al., in \cite{sharma1997eccentric}, is a distance–cum–degree based descriptor defined as $$\displaystyle{\xi^C(G) = \sum_{u \in V}  d_u \;ecc(u)}$$ where $d_u$ is the degree of the vertex $u$ of and $ecc (u)$ is the eccentricity of the vertex $u$ of the graph $G=(V,E)$.

The eccentric connectivity index has been extensively studied in undirected graphs by various authors, to mention a few see~ \cite{morgan2011eccentric,zhou2010eccentric,hauweele2019maximum,zhang2012minimal,zhang2014maximal} and in composite graphs \cite{azari2022study,dovslic2014eccentric} due to its strong correlation with molecular properties.
However, many real-world systems, such as communication, transportation, and biological networks, are inherently directional, leading to asymmetric distances. For this reason, several distance–based indices \cite{knor2016digraphs} and degree–based indices \cite{monsalve2021vertex}  have already been extended to digraphs. As many authors tried to study the interesting concepts in simple graphs to digraphs, it will be interesting to study the concept of eccentric connectivity index in digraphs, which-unlike other topological indices-use both the degree of a vertex and the distances.  Extending this index to directed graphs with respect to the metric maximum distance will be worth studying.

The present work introduces and studies the eccentric connectivity index for strongly connected directed graphs with respect to metric maximum distance defined in \cite{chartrand1997distance}, establishing its basic properties and studying the extremal values of eccentric connectivity index over orientations of the complete graph $K_n$  and related bounds. The article is organised as follows. In section \ref{sec:prelim}, all the basic definitions are provided. In section \ref{sec:ecid}, Eccentric connectivity index of digraph is introduced. In sub-section \ref{subsec:orientedKn}, Eccentric connectivity index of orientations of complete graph $K_n$ is studied. In section \ref{sec:mec}, an attempt is made to study  the maximum eccentric connectivity index of digraph. 
\section{Preliminaries}\label{sec:prelim}
A directed graph (or digraph) $D = (V, A)$ consists of a non-empty finite set $V$ of vertices, and a finite set $A$ of ordered pairs of distinct vertices, called \textit{arcs}, see \cite{bang2008digraphs}. An arc from vertex $u$ to vertex $v$ is denoted by $(u,v)$, indicating a directed edge from $u$ to $v$, and $u$ is called the tail and $v$ is called the head. The \emph{out-degree} of a vertex $u$, denoted by $d^+_u$, is the number of arcs originating from $u$, That is, arcs of the form $(u,v)$. Similarly, the \emph{in-degree} of $u$, denoted by $d^-_u$, is the number of arcs directed towards $u$, That is, arcs of the form $(v,u)$. A digraph $D$ is $r$ regular if and only if $d^+_u=r=d^-_u,$ for all $u \in V(D)$.
 
A \textit{directed walk} (or diwalk) in a digraph $D$ is an alternating sequence $ W = x_1 a_1 x_2 a_2 \ldots x_{k-1} a_{k-1} x_k,
$ where $x_i$ are vertices and $a_j$ are arcs such that each arc $ a_i$ has tail $ x_i $ and head $ x_{i+1}$ for all $ i \in 1,\ldots, k-1$. We denote this by $x_1 \rightarrow{x_2}\rightarrow \cdots \rightarrow{x_{k-1} \rightarrow{x_k}}$. If all the vertices in the walk are distinct, $W$  is called a \textit{directed path} (or dipath). If $x_1 = x_k $, the vertices $ x_1, \ldots, x_{k-1} $ are distinct, and $ k\geq 3 $, then $ W $ is referred to as a \textit{directed cycle} (or dicycle). For a detailed study see \cite{bang2018classes}.  

A digraph $D$ is said to be \textit{strongly connected} (or strong) if, for every pair of distinct vertices $(u, v)$, there exists a directed path from $u$ to $v$. The length of a directed path is the number of arcs it contains.  For any pair of vertices $(u, v)$ in a strongly connected digraph $D$, the shortest directed $u, v$- path is called a directed geodesic, and its length, that is, the number of arcs in it is termed the directed distance, denoted by $\overrightarrow{d}(u, v)$.

A digraph $D = (V,A)$ is symmetric if $(x,y) \in A $ implies $ (y,x) \in A $. An orientation of an undirected digraph $G$ is an oriented graph $D$ obtained from G by replacing every edge $(x,y)$ by either arc $(x,y)$ or arc $(y,x)$.
For an undirected graph $G$ , the complete biorientation of $G$ , denoted by $\overleftrightarrow{G}$ (or $\widehat{G}$) is the symmetric digraph obtained from  $G $ by replacing each edge $ (x,y)$ with the pair of arcs $ (x,y)$  and $ (y,x) $.  The reverse of a strongly connected digraph $D$ obtained by reversing the orientation of all arcs of a digraph $D$ is denoted by $D^-$.

 Similar to out-degree and in-degree of a vertex $u \in D $ we have \textit{out-eccentricity} $ecc^+(u)$ and \textit{in-eccentricity} $ecc^-(u)$.  The out-eccentricity $ecc^+ (u)$ of a vertex $u$ in $D$ is the maximum distance from $u$ to any other vertex of $D$ where as the in-eccentricity $ecc^- (u)$ of a vertex $u$ in $D$ is the maximum distance from any other vertex of $D$ to the vertex $u$. 
 The eccentricity $ecc (u)$ of a vertex $u \in D$ is the maximum of $ecc^+(u)$ and $ecc^-(u)$. That is $ecc(u)= \max  \{ecc^+(u),ecc^-(u)\} , u\in D$, Thus the eccentricity $ecc (u)$ of a vertex $u \in D$ is the maximum distance between $u$ and any other vertex of $D$.
 
 If the digraph is not strongly connected the eccentricity $ecc (u)$ of a vertex $u \in D$ may not be finite, hence we will consider strongly connected digraphs and take the distance as metric ‘\textit{maximum distance}’, abbreviated as \textit{md} was introduced by Chartrand and Tian in \cite{chartrand1997distance}. It is also known as m-distance and is defined as $md(u,v)=\max\{\vec{d}(u,v),\vec{d}(v,u)\}$.

The following definitions are from \cite{chartrand1997distance}.
The $mecc(u)=\max \{md(u,v)\mid v \in V(D)\}$. The $m-radius$, $ mrad(D)$, of a digraph $D$ is defined by  $mrad(D)=\min \{mecc(v)\mid v \in V(D)\}$. The $m-diameter$, $mdiam(D)$, of a digraph $D$ is defined by  $mdiam(D)=\max \{mecc(v)\mid v \in V(D)\}$.
A digraph is self centered if its radius and diameter are equal.
The out-neighborhood $N_D^+(v) = \{u \in V : (v,u) \in A\}$, in-neighborhood $N_D^-(v) = \{w \in V : (w,v) \in A\}$.

\section{Eccentric connectivity index of digraphs}\label{sec:ecid}
In this section we are extending \emph{eccentric connectivity} index to digraphs without loops and parallel arcs.
\begin{definition}
Let $D=(V,A)$ be a digraph, the \emph{eccentric connectivity index} of $D$ is defined as 
$$\displaystyle{\xi^C(D) = \frac{1}{2} \sum_{u \in V} ~~  (d_u^++d_u^-)\;mecc(u) }$$
\end{definition}
To make the idea clear, consider the Example \ref{example:eec}.
\begin{example}\label{example:eec}
Consider the digraph $D$ in Figure \ref{fig:AneqIII}.
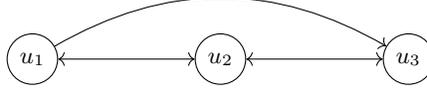
\begin{figure}[H]
\centering
\begin{tikzpicture}[node distance=2.5cm, every node/.style={circle, draw}]
  \node (1) {$u_1$};
  \node (2) [right of=1] {$u_2$};
  \node (3) [right of=2] {$u_3$};

  \draw[<->] (1) -- (2);
  \draw[<->] (2) -- (3);
  \draw[->, out=30, in=150, looseness=1] (1) to (3);
 
\end{tikzpicture}
\caption{The digraph $D$.}
\label{fig:AneqIII}
  \end{figure}
$mecc(u_1)=2$, $mecc(u_2)=1$, $mecc(u_3)=2$, 
$d^+_{u_1}=2$, $d^-_{u_1}=1$, $d^+_{u_2}=2$, $d^-_{u_2}=2$, $d^+_{u_3}=1$, $d^-_{u_3}=2$, and so $\xi^C(D) =8$.
\end{example}
In \cite{monsalve2021vertex}, the authors proved the following theorem in the case of symmetric digraphs.
\begin{theorem}[Theorem 2.4 of \cite{monsalve2021vertex}]
    Let $\phi$ be a symmetric Vetex Degree Based topological index of $G$ and $\widehat{G}$ is the symmetric digraph of $G$. Then $\phi (G) = \phi{(\widehat{G})}$. 
\end{theorem}
Similar result can be proved in the case of a graph and its symmetric digraph also. See the Thm.~\ref{thm:sameecc}.
\begin{theorem}\label{thm:sameecc}
If $G$ be any connected graph, then the eccentric connectivity index of $G$ and its symmetric digraph $\widehat{G}$ are same. That is $\xi^C(G)=\xi^C(\widehat{G})$.
\end{theorem}
\begin{proof}

Let $G=(V,E)$ be a graph and $\widehat{G}=(V,A)$ be its symmetric digraph.\\
Since an edge in $G$ is replaced by a pair of symmetric arcs in $\widehat{G}$, The out-degree as well as in-degree  of each vertex $u \in $ $\widehat{G}$ is same as the degree of the vertex that in $D$.\\
Hence sum of out-degree and in-degree of $u \in V(\widehat{G})$ is double of the degree of  $ u \in V(G)$, that is $d_u^++d_u^-=2d_u$. 
Since $md(u,v)=\max\{\vec{d}(u,v),\vec{d}(v,u)\}=d(u,v)$, we have $mecc(u)=ecc(u)$.  Therefore,
 \begin{align*}
 \xi^C(\widehat{G})&=\frac{1}{2} \sum_{u \in V(\widehat{G})}   (d_u^++d_u^-)mecc(u)\\ 
 &=\frac{1}{2} \sum_{u \in V(G)}   2d_uecc(u) \\ 
 &= \sum_{u \in V(G)}  d_uecc(u) \\
 &=\xi^C(G). 
 \end{align*}
 \qed
\end{proof}

\begin{lemma}\label{lem:scd}
    Let $D=(V,A)$ be a strongly connected digraph and $D^-=(V,A^-)$ be its reverse then $ecc^+(u), u\in V(D)=ecc^-(u), u\in V(D^-)$ and $ecc^-(u), u\in V(D)=ecc^+(u), u\in V(D^-)$.
\end{lemma}
\begin{proof}
    Let $D=(V,A)$ be a strongly connected digraph and let $D^-=(V,A^-)$ be the reverse of $D$ obtained by reversing the orientation of all arcs of $D$.
    Then the out-eccentricity $ecc^+$ of the vertex $u \in V(D)$ is same as  the in-eccentricity $ecc^-$ of the vertex $u \in V(D^-)$ and the in-eccentricity $ecc^+$ of the vertex $u \in V(D)$ is same as the out-eccentricity $ecc^-$ of the vertex $u \in V(D^-)$.  So, for $u \in V(D)$,
    \begin{align*}
       mecc(u) &=\max_{u\in V(D)} \{ecc^+ (u),ecc^- (u)\}\\ 
       &=  \max_{u\in V(D^-)} \{ecc^-(u),ecc^+ (u)\}\\ 
       &= mecc(u)  \text{ for } u \in V(D^-).
  \end{align*}
  \qed
\end{proof}
\begin{theorem}\label{reverse}
    Let $D$ be a strongly connected digraph and let $D^-$ be its reverse. Then $\xi^C(D)=\xi^C(D^-)$.
\end{theorem}
\begin{proof}
    Let $D$ be a strongly connected digraph and let $D^-$ be the reverse of $D$ obtained by reversing the orientation of all arcs of $D$.
    The out-degree of the vertex $u \in V(D)$ becomes the in-degree of the vertex $u \in V(D^-)$ and the in-degree of the vertex $u \in V(D)$ becomes the out-degree of the vertex $u \in V(D^-)$.
    There for the sum of out-degree and in-degree of the vertex $u\in V(D)$ equals to the sum of out-degree and in-degree of the vertex $u\in V(D^-)$.
    That is $d_u^++d_u^-$ is same for $u\in V(D) $ and $u \in V(D^-)$.  Hence by~Lemma~\ref{lem:scd}, 
   \begin{align*}
 \xi^C({D})&=\frac{1}{2} \sum_{u \in V(D)} (d_u^++d_u^-)mecc(u)\\
 &=\frac{1}{2} \sum_{u \in V(D^-)}(d_u^++d_u^-)mecc(u)\\ 
 &= \xi^C({D^-}).
 \end{align*}
\qed
\end{proof}   

\begin{theorem}\label{ecc-connectivity-bound}
If $D$ is a strongly connected digraph with $n$ vertices and $a$ arcs, then $a. mrad(D) \le \xi^C (D) \le a. mdiam(D)$.
\end{theorem}
\begin{proof}
Note that $\displaystyle \sum_{u \in V} \left(d_u^++d_u^-\right) = 2a$.  Also by definition $mrad(D) \le mecc (u) \le mdiam(D) $ for every vertex $u \in V$.  Hence 
\begin{align*}
\xi^C(D)&= \frac{1}{2} \sum_{u \in V}  (d_u^++d_u^-)mecc(u)\\ & \le \frac{1}{2} mdiam(D)\sum_{u \in V}(d_u^++d_u^-)\\ & = \frac{1}{2} mdiam(D)2a =a.mdiam(D)
\end{align*}
\begin{align*}
\xi^C(D)&= \frac{1}{2} \sum_{u \in V}  (d_u^++d_u^-)mecc(u)\\ & \ge \frac{1}{2} mrad(D)\sum_{u \in V}(d_u^++d_u^-)\\ & = \frac{1}{2} mrad(D).2a =a.mrad(D)
\end{align*}
Therefore $a.mrad(D) \le \xi^C(D)\le a.mdiam(D)$.\qed
\end{proof}
Since, for every self centered strongly connected digraph $D$ the radius and diameter of $D$ are the same, it can be seen that the inequality given in the Thm.~\ref{ecc-connectivity-bound} will change to equalty.  We state it as a corrollary of the Thm.~\ref{ecc-connectivity-bound}.
\begin{corollary}Let $D$ be a strongly digraph with $n$ vertices and $a$ arcs. Then  $a. mrad(D) = \xi^C (D) = a. mdiam(D)$ if and only if $D$ is self-centered. 
\end{corollary}
\begin{corollary}\label{complement}
    Let $D$ be a strongly connected digraph with $n \geq 4$ vertices for which the complement is also strongly  connected. Then $\xi^{C}(D) + \xi^{C}(\overline{D}) \geq 2n(n-1)$ with equality if and only if both $D$ and $\overline{D}$ are self-centered digraphs with radius two.
\end{corollary}
\begin{proof}
Let $D$ be a strongly connected digraph on $n$ vertices and $a$ arcs and $\overline{D}$ be its compliment with $\overline{a}$ arcs, which is also strongly connected. 
So $a+\overline {a} =n(n-1)$.  Since both are $D$ and $\overline{D}$ strongly connected, if a vertex $u$ of $D$ has arcs $(u,v)$ and $(v,u)$ for all $v \in D$ then the vertex $u$ will be isolated on the compliment digraph $\overline{D}$, which is not possible.  Hence every vertex of $D$ and $\overline{D}$ has out-degree and in-degree lies between $2$ and $n-2$.  Thus mecc is atleast 2 for every vertex of $D$ and $\overline{D}$.  Therefore  $mrad\ge 2$ for both the digraphs $D$ and   $\overline{D}$.  Then by Thm. \ref{ecc-connectivity-bound}, $\xi^C(D) \ge 2a $ and $\xi^C(\overline{D})  \ge 2\overline{a}$.  Thus 
$\xi^C(D)+\xi^C(\overline{D})  \ge 2a+2\overline{a}=2(a+\overline{a}) \ge 2n(n-1)$.
\qed
\end{proof}
\begin{theorem}
    If $D$ is a $r$-regular digraph with $n$ vertices then $2nr\le \xi^C(D)\le nr(n-r)$. 
\end{theorem}
\begin{proof}
Let $D$ be a $r ( r \ne n-1)$ regular digraph with $n$ vertices namely $u_1,u_2,\ldots,u_n$ then out-degree and in-degree of each  $u_i$ is $r$ and number of arcs in $D$ is $nr$. Hence $d_{u_i}^++d_{u_i}^-=r+r=2r$.  Let out-neighbourhood of the vertex $u_i$ for $i=1,2,\ldots, n$ be $N^+(u_i)=(u_1^i,u_2^i,\ldots,u_r^i)$ where $u_j^i,\; j=1,2,\ldots, r$  are $r$ distinct vertices of $D$.  So, for $k=1,2,\ldots,n$; it can be seen that $\vec{d}(u_i,v_k)=1$ if $v_k\in N^+(u_i)$. If not, then if $v_k\in N^+(u_j^i)$ then  $\vec{d}(u_i,v_k)=2$ , if not, continuing like this we get $\vec{d}(u_i,v_k)\le n-r$ since the digraph $D$ is $r$ regular.  Similarly we get $\vec{d}(v_k,u_i)\le n-r$.
    Thus $md(u_i,u_j)\le n-r$, hence $mecc(u_i)\le n-r$.
    Therefore $mdiam \le n-r$.\\
Since $D$ is $r (r \ne n-1)$ regular no vertex of $D$ have $mecc=1$, so $mrad(D) \ge2 $.\\
  Then by Thm. \ref{ecc-connectivity-bound} we have $2nr \le \xi^C({D})\le nr(n-r)$. \\
  If $r=n-1$, then $D$ is $\overleftrightarrow{K}_n$, hence  $\xi^C({D})=n(n-1)= nr(n-r)$.\\ Equality holds if and only if the digraph is self centered. 
Hence the theorem.\qed
\end{proof}

   \subsection{Eccentric connectivity index of orientations of $K_n$}\label{subsec:orientedKn}
   In this section we study the eccentric connectivity index of orientations of complete graph $K_n$.
Consider the below example to see  different orientations of $K_n$ may have different eccentric connectivity index. We are interested to find which orientation of $K_n$ has maximum and minimum eccentric connectivity index.
 \begin{example}
 Two non isomorphic orientations $T_1$ and $T_2$ of $K_5$ are given in the Fig. \ref{fig:II} 
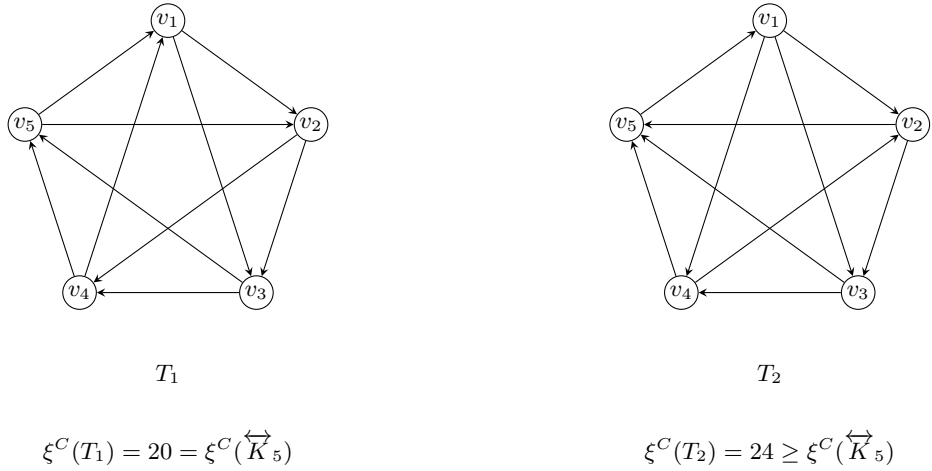
\begin{figure}[H]
\centering
\begin{tikzpicture}[scale=1, >=stealth]

\begin{scope}[xshift=-4cm]
  \foreach \i/\name in {0/v_1,1/v_2,2/v_3,3/v_4,4/v_5} {
    \node[circle, draw, inner sep=1pt, minimum size=4mm] (L\i) at (90-72*\i:2cm) {$\name$};
  }

  \draw[->] (L0) -- (L1);
  \draw[->] (L0) -- (L2);
  \draw[->] (L1) -- (L2);
  \draw[->] (L1) -- (L3);
  \draw[->] (L2) -- (L3);
  \draw[->] (L2) -- (L4);
  \draw[->] (L3) -- (L4);
  \draw[->] (L3) -- (L0);
  \draw[->] (L4) -- (L0);
  \draw[->] (L4) -- (L1);

  \node at (0,-2.7) {$T_1$};
  \node at (0,-3.7) {$\xi^C({T_1})=20=\xi^C(\overleftrightarrow{K}_5) $};

\end{scope}

\begin{scope}[xshift=4cm]
  \foreach \i/\name in {0/v_1,1/v_2,2/v_3,3/v_4,4/v_5} {
    \node[circle, draw, inner sep=1pt, minimum size=4mm] (R\i) at (90-72*\i:2cm) {$\name$};
  }

  \draw[->] (R0) -- (R1);
  \draw[->] (R0) -- (R2);
  \draw[->] (R0) -- (R3);
  \draw[->] (R4) -- (R0);

  \draw[->] (R1) -- (R2);
  \draw[->] (R1) -- (R4);
  \draw[->] (R3) -- (R1);

  \draw[->] (R2) -- (R3);
  \draw[->] (R2) -- (R4);
  \draw[->] (R3) -- (R4);

  \node at (0,-2.7) {$T_2$};
  \node at (0,-3.7) {$\xi^C({T_2})=24 \ge \xi^C(\overleftrightarrow{K}_5{})$};
\end{scope}

\end{tikzpicture}
\caption{Orientations of $K_5$.}
\label{fig:II}
  \end{figure}
 \end{example}
\begin{lemma}
    $\xi^C(\overleftrightarrow{K}_n)=\xi^C(K_n)\le \xi^C(\overrightarrow{K}_n)$, for any strong orientation of $K_n$.
\end{lemma}
\begin{proof}
Let $\overrightarrow{K}_n$ be a strong orientaion of $K_n$, and $\overleftarrow{K}_n$ be its complement.
 By Cor. \ref{complement}, $\xi^{C}(\overrightarrow{K}_n) + \xi^{C}(\overleftarrow{K}_n) \geq 2n(n-1)$.  But, by  Thm. \ref{reverse}, $\xi^{C}(\overrightarrow{K}_n) = \xi^{C}(\overleftarrow{K}_n)$.  Therefore  $2\ \xi^{C}(\overrightarrow{K}_n)  \geq 2n(n-1)$.
 Hence $\xi^{C}(\overrightarrow{K}_n) \ge n(n-1)= \xi^C(K_n)$.
 \qed
\end{proof}

\begin{figure}[H]
\centering

\begin{tikzpicture}[scale=1, >=stealth,
  every node/.style={circle, draw, inner sep=1pt, minimum size=5mm}]

\begin{scope}[xshift=-4cm]
\foreach \i/\lab in {
  0/{$v_1$},1/{$v_2$},2/{$v_3$},3/{$v_4$},4/{$v_5$},5/{$v_6$},6/{$v_7$}
} {
  \node (A\i) at (90-360/7*\i:2cm) {\lab};
 
}

\foreach \i in {0,...,6} {
  \pgfmathtruncatemacro{\j}{mod(\i+1,7)}
  \draw[->,  thick, red] (A\i) -- (A\j);
}
\foreach \i in {0,...,6} {
  \pgfmathtruncatemacro{\j}{mod(\i+2,7)}
  \draw[->,  thick, brown] (A\i) -- (A\j);
}
\foreach \i in {0,...,6} {
  \pgfmathtruncatemacro{\j}{mod(\i+3,7)}
  \draw[->,  thick, blue] (A\i) -- (A\j);
}

\end{scope}

\begin{scope}[xshift=4cm]
\foreach \i/\lab in {
  0/{$v_1$},1/{$v_2$},2/{$v_3$},3/{$v_4$},4/{$v_5$},5/{$v_6$},6/{$v_7$}
} {
  \node (B\i) at (90-360/7*\i:2cm) {\lab};
 
}
\node[fill=yellow!20] (B7) at (0,0) {$v_8$};

\foreach \i in {0,...,6} {
  \pgfmathtruncatemacro{\j}{mod(\i+1,7)}
  \draw[->, thick, red] (B\i) -- (B\j);
}
\foreach \i in {0,...,6} {
  \pgfmathtruncatemacro{\j}{mod(\i+2,7)}
  \draw[->, thick, brown] (B\i) -- (B\j);
}
\foreach \i in {0,...,6} {
  \pgfmathtruncatemacro{\j}{mod(\i+3,7)}
  \draw[->, thick, blue] (B\i) -- (B\j);
}

\foreach \a in {0,2,4,6} {
  \draw[->, very thick, black] (B7) -- (B\a);
}
\foreach \b in {1,3,5} {
  \draw[->, very thick, black] (B\b) -- (B7);
}

\end{scope}

\end{tikzpicture}
\caption{$\overrightarrow{K}_n $ for $n=7$ and $n=8$.}
\label{fig:III}
  \end{figure}
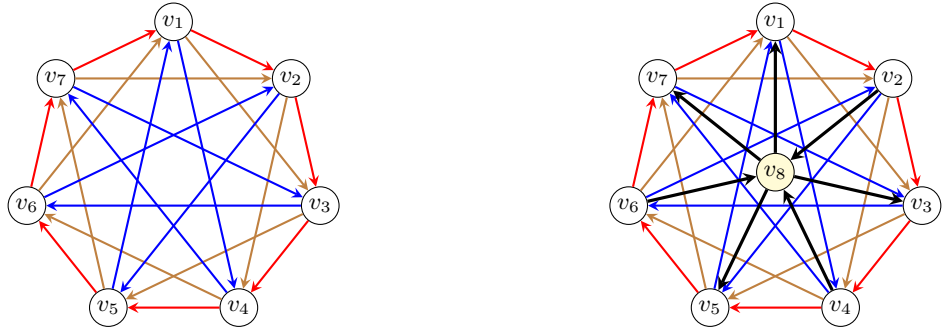

\begin{theorem}\label{thm:Knorient}
For $n\ge3$ there exist an orientation $\overrightarrow{K}_n$ of the complete graph $K_n$ such that  $\xi^C(\overrightarrow{K}_n)=\xi^C(K_n)=\xi^C(\overrightarrow{C}_n)$.
\end{theorem}
 \begin{proof}
 To prove the theorem, we give different orientations for $K_n$ when $n$ is odd and when $n$ is even.  Let $v_1,v_2,\ldots,v_n$ be the vertices of $K_n$. \\
 \noindent
 \textbf{Case 1:} $n$ is odd.\\
Construct the orientation $\overrightarrow{K}_n$ of the complete graph $K_n$ as follows.
\begin{description}
    \item[Step 1:] The cycle $v_iv_{i+1} \ldots v_n v_1\ldots v_{i-1} v_i$ in $K_n$ is replaced by the directed cycle $v_i\rightarrow{v_{i+1}}\rightarrow\cdots\rightarrow v_n\rightarrow v_1\rightarrow \cdots \rightarrow v_{i-1} \rightarrow v_i $.
    \item[Step 2:] The cycle $v_i v_{i+2} \ldots v_{n-1} v_1 \ldots v_{i-2} v_i$ in $K_n$ is replaced by the directed cycle $v_i\rightarrow{v_{i+2}}\rightarrow\cdots\rightarrow v_{n-2}\rightarrow v_1\rightarrow \cdots \rightarrow v_{i-2} \rightarrow v_i $.
\end{description}
Continuing like this the process will terminate after $\frac{n-1}{2}$ steps with all the vertices have out-degree and in-degree  $\frac{n-1}{2}$.
 So, the sum of out-degree and in-degree of each vertex is $n-1$. By construction every vertex $v_j, j \ne i$ is in $N_D^+(v_i)$ or in $N_D^{++}(v_i)$.  Hence $mecc$ of each vertex is $2$.
 Therefore, $\xi^C(\overrightarrow{K}_n)= n(n-1)=\xi^C(K_n)=\xi^C(\overrightarrow{C}_n)$.

 \noindent   
\textbf{Case $2$:} $n$ is even.\\
Consider $\overrightarrow{K}_{n-1}$  with vertices $v_1,v_2,\ldots, v_{n-1}$ and add arcs from $v_n$ to $v_i$ for $i=1,3,\ldots, n-1$ and $v_i$ to $v_n$ for $i=2,4,\ldots,n$.  By the construction it can be seen that, adding $v_n$ will not change the mecc of $v_1,v_2,\ldots v_{n-1}$ 
so that $mecc(v_i)=2$ for all $i=1,2,\ldots,n-1$. \\
If $i$ is odd, then $\overrightarrow{d}(v_n,v_i)=1$ for all $i=1,3,\ldots,n-1$ and since $v_i \rightarrow{v_{i+1}}$ and $v_{i+1} \rightarrow{v_{n}}$ for $i\ne n-1$ and Since $(v_{n-1},v_2), (v_2,v_n)$ are arcs, we have $\overrightarrow{d}(v_{n-1},v_n)=2$.\\ 
Now for $i$ even, $\overrightarrow{d}(v_i,v_n)=1$ for all $i=2,4,\ldots,n$ and  since $v_n \rightarrow{v_{i-1}}$ and $v_{i-1} \rightarrow{v_{i}}$ we have $\overrightarrow{d}(v_n, v_i)=2$.\\
Thus $mecc(v_i)=2$ for all vertices $v_i$.  Therefore, $\xi^C(\overrightarrow{K}_n)=n(n-1)=\xi^C(K_n)=\xi^C(\overrightarrow{C}_n)$.\qed

 \end{proof}
 To explain the Thm.~\ref{thm:Knorient} see the Fig.~\ref{fig:III}, where $\overrightarrow{K}_8$ is obtained from  $\overrightarrow{K}_7$ by adding arcs from $v_8$ and to $v_8$ in $\overrightarrow{K}_8$ as described in Case 2.

\begin{theorem}\label{star}
    If $D=(V,A)$ be a strongly connected digraph of order $n, n \ge 4$ then $\xi^C(D)\ge 3(n-1)$, and equality holds if and only if $D$ is the biorientation of star graph $\overleftrightarrow{S}_n$. 
\end{theorem}
\begin{proof}
Let $D=(V,A)$ be a strongly connected digraph of order $n, n \ge 4$ and  $A=\{ u \in V \mid d^+_u=d^-_u=n-1 \}$ and $B=V-A$. \\
Let $ \mid A \mid =k$ then $ \mid B \mid =n-k$.\\
Since every vertex of $A$ has out-degree and in-degree  $n-1$, $mecc$ of each vertex in $A$  is one and $mecc(u) \ge 2$ for each $u \in B$.\\ 
Let $d_u = d_u^+ + d_u^-$.\\
We have $$ \displaystyle{\sum_{u \in V}  d_u=2m  \implies \sum_{u \in A}  d_u + \sum_{u \in B}  d_u=2m}$$\\ 
\begin{equation}{\label{1}}
 \displaystyle{\implies \sum_{u \in B}d_u =2m-  \sum_{u \in A}  d_u }
 \end{equation}
Now $$\displaystyle{\sum_{u \in A}  d_u = \sum_{u \in A}2(n-1)}=2k(n-1).$$ 
Then Eqn.(\ref{1}) becomes $$\displaystyle{\sum_{u \in B}  d_u =2m-2k(n-1)}$$
Now,
\begin{align*}
\xi^C(D)&= \frac{1}{2} \left [ \sum_{u \in A} d_u mecc(u)+\sum_{u \in B} d_u mecc(u)\right ]\\ & \ge \frac{1}{2} \left [ \sum_{u \in A} 2(n-1).1 +\sum_{u \in B} d_u.2\right ] \\  & = \left [  k(n-1) +(2m-2k(n-1))\right ]  = 2m-k(n-1).
\end{align*}
That is, \begin{equation}{\label{2}}
  \xi^C(D)\ge 2m-k(n-1). 
  \end{equation}
\textbf{Case 1.} $k=0 \:(A=\phi)$.\\
Then from Eqn. (\ref{2}), $\xi^C(D)\ge 2m\ge 3(n-1)$ if $m\ge \frac{3}{2}(n-1)$.\\
If $m< \frac{3}{2}(n-1)$ partition the set $B$ into $B_1=\{ u \in V \mid d^+_u=d^-_u=1 \}$ where $ \mid B_1 \mid =\ell$ and $B_2=B-B_1$ then $ \mid B_2 \mid =n-\ell$.
Now $d^+_u+d^-_u \ge 3$ and $mecc(u)\ge 2$ for every vertex $u \in B_2$. Therefore,
\begin{align*}
\xi^C(D)&= \frac{1}{2} \left [ \sum_{u \in B_1} d_u mecc(u)+\sum_{u \in B_2} d_u mecc(u)\right ]\\ 
& \ge \frac{1}{2} \left [ \sum_{u \in B_1} 2.2 +\sum_{u \in B_2} 3.2\right ] \\  
& = 3n-\ell  \ge 3n -3 \text{ for } \ell\leq 3.
\end{align*} 
Now, for $\ell > 3$ we can see that $mecc(u)\ge 3$ for every $u\in B_2$.
Hence 
\begin{align*}
\xi^C(D)&= \frac{1}{2} \left [ \sum_{u \in B_1} d_u mecc(u)+\sum_{u \in B_2} d_u mecc(u)\right ]\\
& \ge \frac{1}{2} \left [ \sum_{u \in B_1} 2.3 +\sum_{u \in B_2} 3.2\right ] \\  
& = 3n  \ge 3n -3.
\end{align*}

\noindent
\textbf{Case 2.} $k\ge 1$.\\
Here all vertices of $B$ has out-degree and in-degree at least $k$ (because $k$ vertices have out-degree$=n-1=$in-degree, so there are $k$ arcs at least from $u$ and to $u$ for every $u \in B$). So we have $2m\ge k.2(n-1)+(n-k)2k$. From Eqn.~(\ref{2}),  $\xi^C(D)\ge k[3n-2k-1]$. \\ Obviously, the function $f(x)=x(3n-2x-1)$ with $1 \leq x \leq n-1$ attains its minimum value for 
$x = 1$ or $x = n-1$.  Note that $f(1)  = 3(n-1) \le n(n-1)=f(n)$.   Hence,  $\xi^{C}(D) \;\geq\; f(1) \;=\; 3(n-1)$.
Thus the equality holds when $k=1$ and $m=2(n-1)$, that is $D=\overleftrightarrow{S}_n$. 
 \qed
\end{proof}
 
\section{Maximum Eccentric connectivity index of digraphs}\label{sec:mec}
  Minimum eccentric connectivity of a digraph on $n$ vertices is attained for the biorientation of the star graph (as in the case of a simple graph) See~Thm. \ref{star}.  In this section, an attempt is made to find the maximum eccentric connectivity of a digraph on $n$ vertices. Many studies where done to find the upperbounds of eccentric connectivity index of graphs. Thm. 9 of \cite{hauweele2019maximum} found maximum eccentric connectivity index of graphs on $n$ vertices.

\begin{figure}[H]
\centering
\begin{tikzpicture}[node distance=1.5cm, 
  every node/.style={circle, draw, minimum size=9mm}, >=stealth]

  \node (1) {$u_1$};
  \node (2) [right of=1] {$u_2$};
  \node (3) [right of=2] {$u_3$};
  \node (4) [right of=3] {$u_4$};
  \node (5) [right of=4] {$\cdots$};
  \node (6) [right of=5] {$u_{n-1}$};
  \node (7) [right of=6] {$u_n$};

  \foreach \i/\j in {1/2,2/3,6/7}{
    \draw[<->] (\i) -- (\j);
  }
  \draw[dotted] (3) -- (4);
  \draw[dotted] (4) -- (5);
  \draw[dotted] (5) -- (6);

  \foreach \j in {3,4,6,7} {
    \draw[->, bend left=30] (1) to (\j);
  }
  \foreach \j in {4,6,7} {
    \draw[->, bend left=30] (2) to (\j);
  }
  \foreach \j in {5,6,7} {
    \draw[->, bend left=25] (3) to (\j);
  }

  \draw[dotted,->, bend left=30] (4) to (6);
  \draw[dotted,->, bend left=30] (4) to (7);
  \draw[dotted,->, bend left=30] (5) to (7);

\end{tikzpicture}
\caption{The biorientation $\overleftrightarrow{P}_n$ along with additional arcs $(u_i,u_j)$ for all $j \ge i+2$.}
\label{fig:AneqIV}
\end{figure}
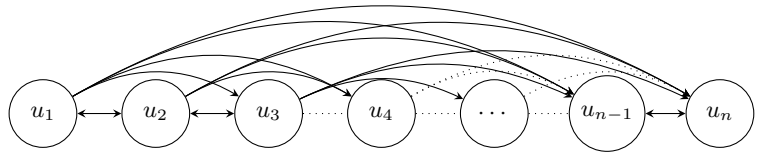

 
If $P_n=u_1u_2\ldots u_n$ is any (simple) path, then $P_n^*$ is the digraph obtained from the biorientation of $P_n$ by adding exactly one of  the arcs $(u_1,u_n)$ or $(u_n,u_1)$.
\begin{theorem}\label{ecc-p*_n}
 If $P_n$ is any simple path, then $\xi^C(P_n)= \xi^C(\overleftrightarrow{P_n}) \le \xi^C(P_n^*)$. 
\end{theorem}
\begin{proof}
Let ${P_n}=u_1u_2\ldots u_n$ be a dipath and $\overleftrightarrow{P}_n$ be its biorientation.
Then by Thm \ref{thm:sameecc}    we have $\xi^C(P_n)= \xi^C(\overleftrightarrow{P_n})$.  Since $P_n^*$ is obtained from the biorientation of the path $P_n$ by adding exactly one arc from the last vertex to the first vertex in any direction (say arc $(u_1,u_n))$,  mecc of every vertex remains the same where as out-degree $u_1$ and in-degree of $u_n$ increases by one.  Therefore  
\begin{align*}
\xi^C(P^*_n)&= \frac{1}{2} \sum_{u \in V(P^*_n)}  (d_u^++d_u^-)mecc(u)\\ 
&= \frac{1}{2} [(d^+_{u_{1}}+d^-_{u_{1}})mecc(u_1)+ \cdots + (d^+_{u_{n}}+d^-_{u_{n}})mecc(u_n)] \\ 
&= \frac{1}{2}  \left [ \sum_{u \in V(P_n)}  (d_u^++d_u^-)mecc(u)+mecc(u_1)+mecc(u_n)\right] \\ 
&= \frac{1}{2} \sum_{u \in V(P_n)}  (d_u^++d_u^-)mecc(u)+ \frac{1}{2}[mecc(u_1)+mecc(u_n)]\\ 
&= \xi^C(P_n)+ \frac{1}{2} [n-1+n-1]= \xi^C(P_n)+n-1
\end{align*}
So, $\xi^C(P^*_n) \ge \xi^C(P_n)=\xi^C(\overleftrightarrow{P}_n)$, and hence the theorem.  
\qed
\end{proof}
The Thm. \ref{ecc-p*_n} imply, one can increase the eccentric connectivity index of a digraph by adding arcs without changing the the eccentricities of the vertices. So the problem here is to find the maximum number of arcs that can be added to $\overleftrightarrow{P}_n$ so that meccentricity remains the same but as number of arcs increases the eccentric connectivity index become maximum.  To see, this consider the Fig. \ref{fig:AneqIV}. 

Here $P^+_n = (V,A)$ where $ V=\{u_1,u_2,\dots,u_n\}, \; A=\{(u_i,u_{i+1}),(u_{i+1},u_i)\mid 1\le i \le n-1\}\cup\{(u_i,u_j)\mid j\ge i+2\}$.  
Since we added arcs $(u_i,u_j)$ for $ j\ge i+2$ (only in one direction)  $\overrightarrow{d}(u_i, u_j)=1$ but $\overrightarrow{d}(u_j, u_i)$ remains the same for $ j\ge i+2$. Hence $md(u_i,u_j)$ will not change.
The out-degree of the vertices $u_1, u_2, \ldots,u_{n-3},u_{n-2}$ increases by $n-2, n-3, \dots,2,1 $ and in-degree of the vertices $u_3, u_4, \ldots,u_{n-1},u_n$ increases by $1,2, \ldots, n-3, n-2$ respectively.
 \begin{align*}
 \begin{split}
     \xi^C(P^+_n)&= \frac{1}{2} \sum_{u \in V(P^+_n)}  (d_u^++d_u^-)mecc(u)\\ 
&= \frac{1}{2} [(d^+_{u_{1}}+d^-_{u_{1}})mecc(u_1)+ \cdots + (d^+_{u_{n}}+d^-_{u_{n}})mecc(u_n)] \\ 
&= \frac{1}{2}  \left [ \sum_{u \in V(P_n)}  (d_u^++d_u^-)mecc(u)+(n-2)mecc(u_1)+(n-3)mecc(u_2)+ \cdots \right.\\&\left.\quad + (n-3)mecc(u_{n-1})+(n-2)mecc(u_n)\right]  
\end{split}
\end{align*}
When $n$ is odd, we have\\
\begin{align*}
\begin{split}
\xi^C(P^+_n) &= \frac{1}{2} \sum_{u \in V(P_n)}  (d_u^++d_u^-)mecc(u)+\\
&\quad\frac{1}{2}\left [ 2(n-2)(n-1)+2(n-3)\left ((n-2)+(n-3)+\cdots \right ) +\frac{n-1}{2}\right]\\
&=  \xi^C(P_n)+ \frac{n^3-4n^2+6n-3}{4}.
\end{split}
\end{align*}

When, $n$  even,

\begin{align*}
\begin{split}
\xi^C(P^+_n) &= \frac{1}{2} \sum_{u \in V(P_n)}  (d_u^++d_u^-)mecc(u)+\\
&\quad \frac{1}{2}\left [ 2(n-2)(n-1)+2(n-3)\left ((n-2)+(n-3)+\cdots+\frac{n}{2} \right )\right]\\ 
&=  \xi^C(P_n)+\frac{3n^3-5n^2-12n+16}{8}.
\end{split}
\end{align*}


\bibliographystyle{plain}
\bibliography{ref}

\begin{thebibliography}{10}

\bibitem{azari2022study}
Mahdieh Azari.
\newblock A study of a new variant of the eccentric connectivity index for
  composite graphs.
\newblock {\em Journal of Discrete Mathematical Sciences and Cryptography},
  25(8):2583--2596, 2022.

\bibitem{bang2018classes}
J{\o}rgen Bang-Jensen and Gregory Gutin.
\newblock {\em Classes of directed graphs}, volume~11.
\newblock Springer, 2018.

\bibitem{bang2008digraphs}
J{\o}rgen Bang-Jensen and Gregory~Z Gutin.
\newblock {\em Digraphs: theory, algorithms and applications}.
\newblock Springer Science \& Business Media, 2008.

\bibitem{chartrand1997distance}
Gary Chartrand and Songlin Tian.
\newblock Distance in digraphs.
\newblock {\em Computers \& Mathematics with Applications}, 34(11):15--23,
  1997.

\bibitem{dovslic2014eccentric}
Tomislav Do{\v{s}}lic and Mahboubeh Saheli.
\newblock Eccentric connectivity index of composite graphs.
\newblock {\em Util. Math}, 95:3--22, 2014.

\bibitem{hauweele2019maximum}
Pierre Hauweele, Alain Hertz, Hadrien M{\'e}lot, Bernard Ries, and Gauvain
  Devillez.
\newblock Maximum eccentric connectivity index for graphs with given diameter.
\newblock {\em Discrete Applied Mathematics}, 268:102--111, 2019.

\bibitem{knor2016digraphs}
Martin Knor, Riste {\v{S}}krekovski, and Aleksandra Tepeh.
\newblock Digraphs with large maximum wiener index.
\newblock {\em Applied Mathematics and Computation}, 284:260--267, 2016.

\bibitem{monsalve2021vertex}
Juan Monsalve and Juan Rada.
\newblock Vertex-degree based topological indices of digraphs.
\newblock {\em Discrete Applied Mathematics}, 295:13--24, 2021.

\bibitem{morgan2011eccentric}
MJ~Morgan, Simon Mukwembi, and Henda~C Swart.
\newblock On the eccentric connectivity index of a graph.
\newblock {\em Discrete Mathematics}, 311(13):1229--1234, 2011.

\bibitem{sharma1997eccentric}
Vikas Sharma, Reena Goswami, and AK~Madan.
\newblock Eccentric connectivity index: A novel highly discriminating
  topological descriptor for structure- property and structure- activity
  studies.
\newblock {\em Journal of chemical information and computer sciences},
  37(2):273--282, 1997.

\bibitem{zhang2014maximal}
Jian-bin Zhang, Zhong-zhu Liu, and Bo~Zhou.
\newblock On the maximal eccentric connectivity indices of graphs.
\newblock {\em Applied Mathematics-A Journal of Chinese Universities},
  29(3):374--378, 2014.

\bibitem{zhang2012minimal}
Jianbin Zhang, Bo~Zhou, and Zhongzhu Liu.
\newblock On the minimal eccentric connectivity indices of graphs.
\newblock {\em Discrete Mathematics}, 312(5):819--829, 2012.

\bibitem{zhou2010eccentric}
Bo~Zhou and Zhibin Du.
\newblock On eccentric connectivity index.
\newblock {\em arXiv preprint arXiv:1007.2235}, 2010.

\end{thebibliography}
\end{document}